\documentclass[12pt]{article}%
\usepackage{amsfonts}
\usepackage{sw20bams}
\usepackage{amsmath}
\usepackage{amssymb}
\usepackage{graphicx}%
\providecommand{\U}[1]{\protect\rule{.1in}{.1in}}
\begin{document}

\title{Four Vignettes on Apparent Size}
\author{Steven Finch}
\date{December 10, 2018}
\maketitle

\begin{abstract}
Problems in optimization and geometric probability are discussed, all
connected with angles subtended at an observer's eye by an object at a
distance. \ Several of these remain unsolved.

\end{abstract}

\footnotetext{Copyright \copyright \ 2018 by Steven R. Finch. All rights
reserved.}A\ well-known exercise in calculus asks how far an observer should
stand from a wall so as to maximize the apparent size of an object (a
painting) on the wall \cite{Do-rail, LS-rail, Ta-rail}. \ Assume that the
floor is represented by the positive $x$-axis and that the object is modeled
by a subinterval $[r-1,r+1]$ of the positive $y$-axis (hence $r>1$). \ The
height of the observer is negligible. \ By the Law of Cosines, the angle
$\Omega(r,x)$ subtended at his/her eye by the line segment is
\[
\Omega=\arccos\left(  \frac{-4+\left[  x^{2}+(r-1)^{2}\right]  +\left[
x^{2}+(r+1)^{2}\right]  }{\sqrt{x^{2}+(r-1)^{2}}\sqrt{x^{2}+(r+1)^{2}}%
}\right)
\]
and this is largest when%
\[%
\begin{array}
[c]{ccc}%
x=x_{\text{max}}=\sqrt{r^{2}-1}, &  & \Omega_{\text{max}}=\arccos\left(
\dfrac{\sqrt{r^{2}-1}}{r}\right)  .
\end{array}
\]
For future reference, we note that $x_{\text{max}}\sim r$ and, for fixed $x$,
\[
\Omega(r,x)\sim\frac{2x}{r^{2}}%
\]
as $r\rightarrow\infty$. \ The preceding serves as inspiration for the
remainder of this paper:\ two vignettes on optimization and two on geometric
probability (Table 1).

\begin{center}%
\begin{tabular}
[c]{ccc}%
{\includegraphics[
height=1.1938in,
width=1.1938in
]%
{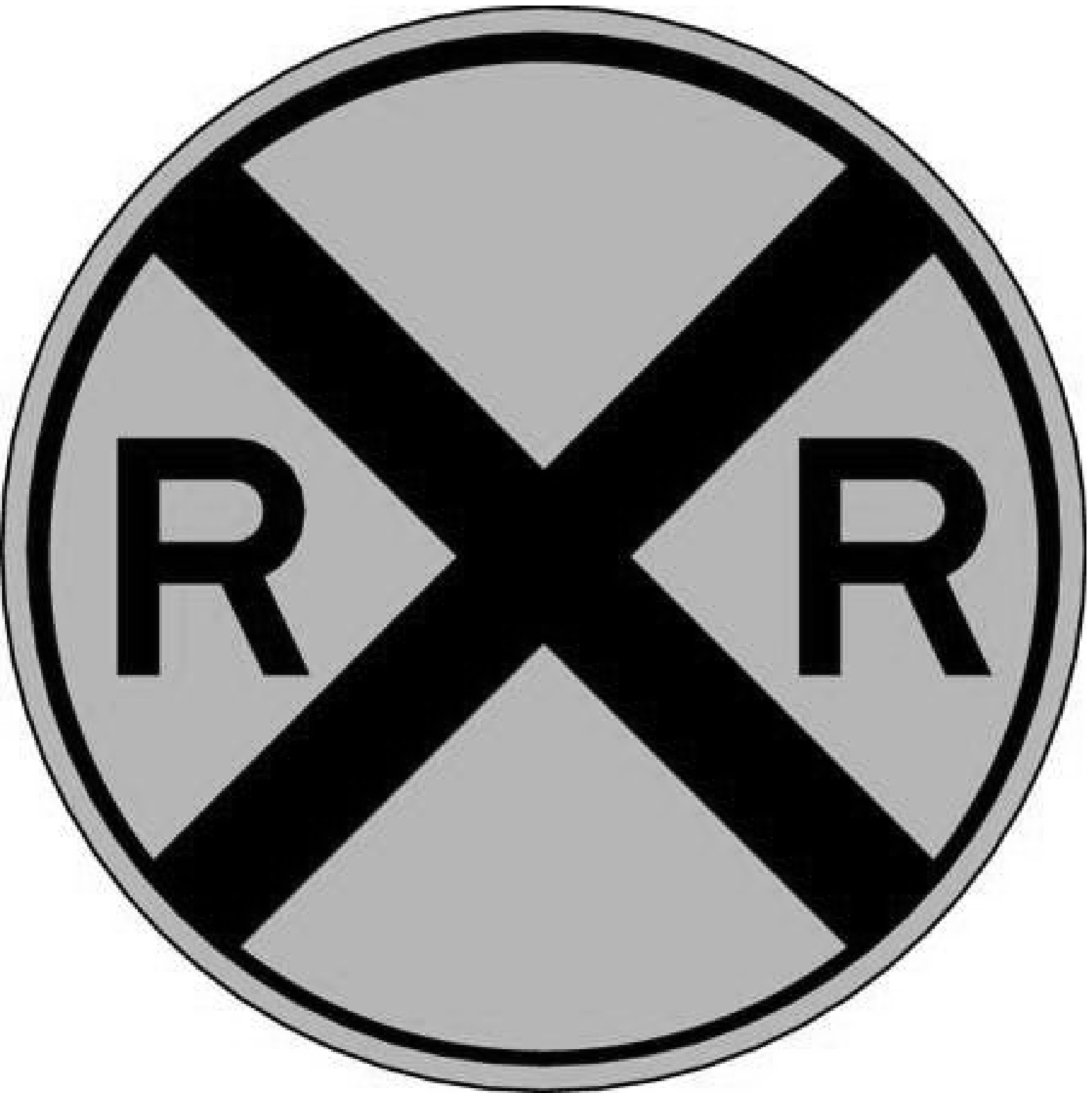}%
}%
&  &
{\includegraphics[
height=1.3317in,
width=1.7509in
]%
{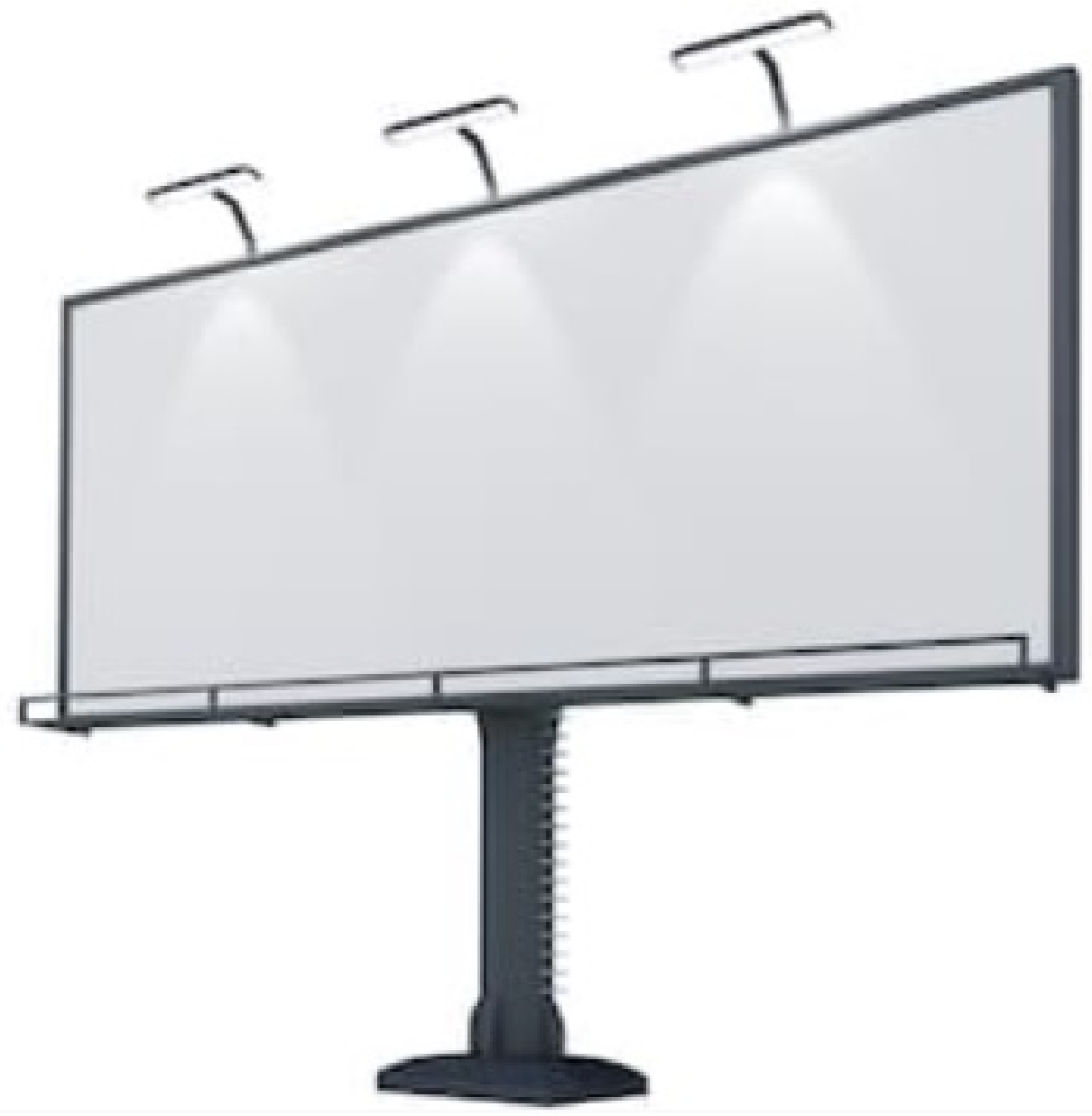}%
}%
\\
&  & \\%
{\includegraphics[
height=1.0984in,
width=2.3943in
]%
{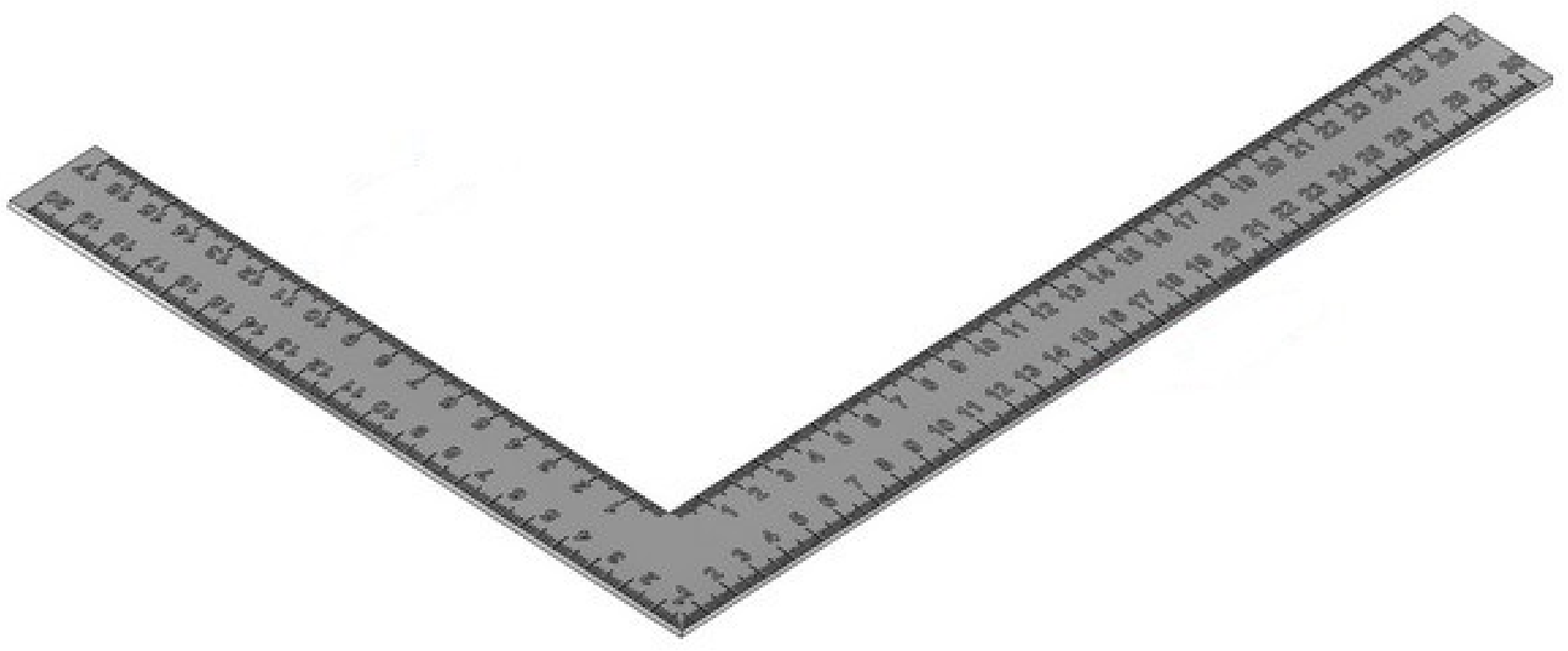}%
}%
&  &
{\includegraphics[
height=1.0843in,
width=1.0843in
]%
{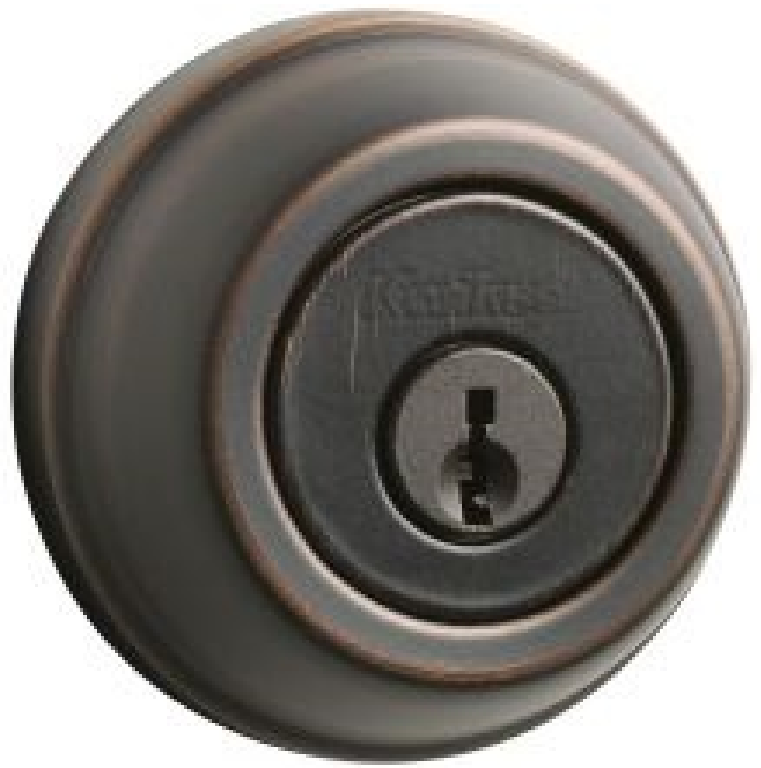}%
}%
\end{tabular}

Table 1: Starting from upper left and moving clockwise:\ a railroad advance
warning sign, a rectangular billboard, a circular deadbolt with keyhole, and a
carpenter's steel framing square.

\end{center}

\section{Railroad Advance Warning Sign}

How far should an observer (a motorist) stand from an object (a circular road
sign) so as to maximize its apparent size? \ Assume that the road is
represented by the positive $x$-axis and that the object is modeled by a disk
of radius $1$ in the $yz$-plane, centered at a fixed distance $r$ from the
origin (again $r>1$). \ The height of the observer is negligible (Figure 1).
\ \textquotedblleft Apparent size\textquotedblright\ is here quantified as the
solid angle $\Omega(r,x)$ subtended by the disk at the observer's eye. \ We
wish to determine the vehicle location $x=x_{\text{max}}$ for which
$\Omega(r,x)$ is maximized, i.e., the sign has greatest visual impact on the
motorist (occupying the largest field of view).%
\begin{figure}[ptb]%
\centering
\includegraphics[
height=3.6148in,
width=4.2026in
]%
{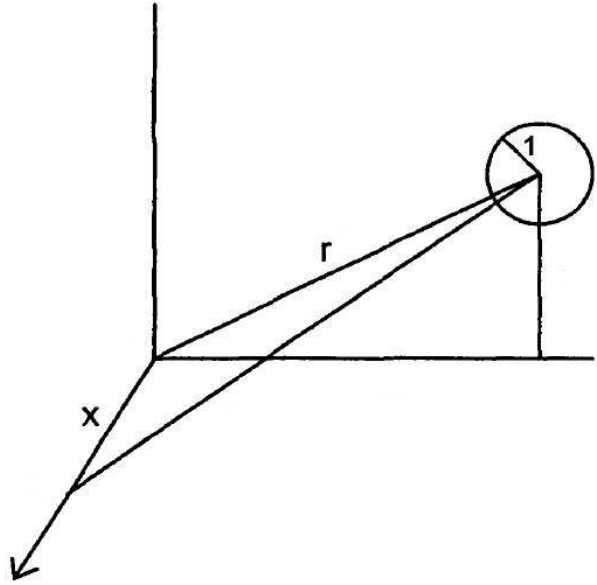}%
\caption{Observer is at $(x,0,0)$; object in $yz$-plane is separated by $r$
from origin.}%
\end{figure}

This problem was solved in \cite{Jo-rail} via a Fourier-Legendre series
approximation for $\Omega$. \ An exact formula, however, was mentioned by
Maxwell \cite{M1-rail} in 1873, explicitly given by Tallqvist \cite{M2-rail}
in 1931, and subsequently rediscovered several times \cite{M3-rail, M4-rail,
M5-rail, M6-rail, M7-rail, M8-rail, M9-rail, M10-rail, M11-rail, M12-rail,
M13-rail, M14-rail}:%

\begin{align*}
\Omega(r,x)  &  =\pi\left(  1+\operatorname*{sgn}(1-r)\right)  -\frac
{2x}{\sqrt{(1+r)^{2}+x^{2}}}\\
&  \cdot\left(  K\left[  \frac{4r}{(1+r)^{2}+x^{2}}\right]  +\frac{1-r}%
{1+r}\Pi\left[  \frac{4r}{(r+1)^{2}},\frac{4r}{(1+r)^{2}+x^{2}}\right]
\right)
\end{align*}
where%
\[
1+\operatorname*{sgn}(1-r)=\left\{
\begin{array}
[c]{lll}%
2 &  & \text{if }r<1,\\
1 &  & \text{if }r=1,\\
0 &  & \text{if }r>1;
\end{array}
\right.
\]%
\[
K[m]=%
{\displaystyle\int\limits_{0}^{1}}
\dfrac{dt}{\sqrt{1-t^{2}}\,\sqrt{1-m\,t^{2}}}\,
\]
is the complete elliptic integral of the first kind; and%
\[
\Pi\left[  n,m\right]  =%
{\displaystyle\int\limits_{0}^{1}}
\dfrac{dt}{\left(  1-n\,t^{2}\right)  \,\sqrt{1-t^{2}}\,\sqrt{1-m\,t^{2}}}%
\]
is the complete elliptic integral of the third kind (Figure 2). \ This formula
enables precise asymptotics that could not be deduced in \cite{Jo-rail}. \ We
have%
\[
x_{\text{max}}\sim\frac{1}{\sqrt{2}}r-\frac{7\sqrt{2}}{24}%
\]
as $r\rightarrow\infty$, which implies that%
\[
x_{\text{max}}\approx(0.707...)r-(0.412...)
\]
for large $r$ and which is an improvement on $x_{\text{max}}\approx
(0.7)r-(0.1)$ given earlier. \ We also have, for fixed $x>0$,
\[
\Omega(r,x)\sim\frac{\pi x}{r^{3}}%
\]
which plays a role in Koopman's \textquotedblleft inverse cube law of
detection\textquotedblright\ from search theory \cite{K1-rail, K2-rail,
K3-rail}.%
\begin{figure}[ptb]%
\centering
\includegraphics[
height=3.4612in,
width=6.1228in
]%
{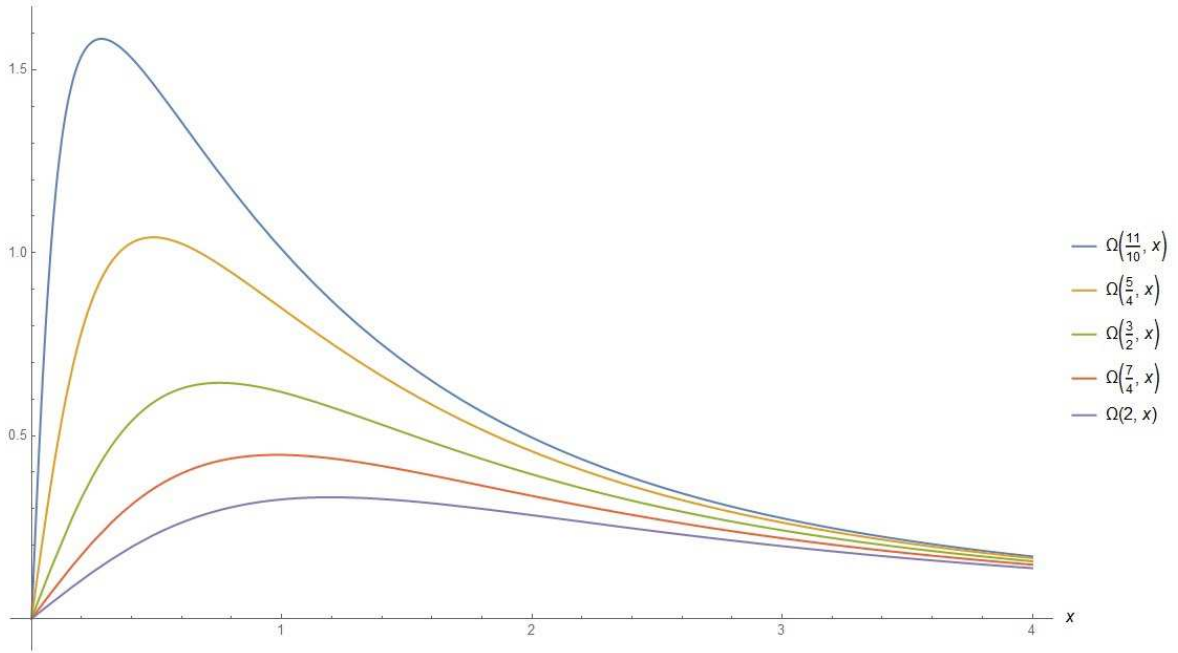}%
\caption{Solid angle as a function of $x$, for certain values of $r$.}%
\end{figure}

\section{Square versus Rectangular Billboards}

The position of the circular sign in the preceding vignette was unspecified:
it sufficed to require only that its center be at distance $r$ from the
origin. \ Here, for a rectangular sign, we must further stipulate that its
center possess $yz$-coordinates%
\[
\left(  \frac{r}{\sqrt{2}},\frac{r}{\sqrt{2}}\right)
\]
that is, it lies on the diagonal line $z=y$. \ Moreover, given that the
rectangle has length $\ell\geq1$ and width $1/\ell$, its vertices are assumed
to be at%
\[%
\begin{array}
[c]{ccc}%
(y_{1},z_{2})=\left(  \dfrac{r}{\sqrt{2}}-\dfrac{\ell}{2},\dfrac{r}{\sqrt{2}%
}+\dfrac{1}{2\ell}\right)  , &  & (y_{2},z_{2})=\left(  \dfrac{r}{\sqrt{2}%
}+\dfrac{\ell}{2},\dfrac{r}{\sqrt{2}}+\dfrac{1}{2\ell}\right)  ,
\end{array}
\]%
\[%
\begin{array}
[c]{ccc}%
(y_{1},z_{1})=\left(  \dfrac{r}{\sqrt{2}}-\dfrac{\ell}{2},\dfrac{r}{\sqrt{2}%
}-\dfrac{1}{2\ell}\right)  , &  & (y_{2},z_{1})=\left(  \dfrac{r}{\sqrt{2}%
}+\dfrac{\ell}{2},\dfrac{r}{\sqrt{2}}-\dfrac{1}{2\ell}\right)  .
\end{array}
\]
Letting the observer be at $(x_{0},y_{0},z_{0})=(x,0,0)$, we have
\cite{C1-rail, C2-rail, C3-rail, C4-rail, C5-rail}%
\begin{align*}
\Omega &  =\arctan\left[  \frac{(y_{2}-y_{0})(z_{2}-z_{0})}{x_{0}\sqrt
{x_{0}^{2}+(y_{2}-y_{0})^{2}+(z_{2}-z_{0})^{2}}}\right]  -\arctan\left[
\frac{(y_{1}-y_{0})(z_{2}-z_{0})}{x_{0}\sqrt{x_{0}^{2}+(y_{1}-y_{0}%
)^{2}+(z_{2}-z_{0})^{2}}}\right] \\
&  -\arctan\left[  \frac{(y_{2}-y_{0})(z_{1}-z_{0})}{x_{0}\sqrt{x_{0}%
^{2}+(y_{2}-y_{0})^{2}+(z_{1}-z_{0})^{2}}}\right]  +\arctan\left[
\frac{(y_{1}-y_{0})(z_{1}-z_{0})}{x_{0}\sqrt{x_{0}^{2}+(y_{1}-y_{0}%
)^{2}+(z_{1}-z_{0})^{2}}}\right]
\end{align*}
and consequently

\begin{itemize}
\item if $x=1$ and $r=1$, then $\ell_{\text{max}}=1$ (a square)

\item if $x=1$ and $r=9/8$, then $\ell_{\text{max}}%
=1.0687058269964383331915489...$

\item if $x=1$ and $r=5/4$, then $\ell_{\text{max}}%
=1.6697745993679378292165263...$
\end{itemize}

\noindent This may be surprising to some readers: the rectangle of largest
solid angle (with fixed observer and fixed center) need not be a square
\cite{Ta-rail, Jo-rail}. \ The asymmetry (due to off-center location of
observer relative to object) is responsible for this phenomena. \ For
$r>1.309987792$, the inequality $r/\sqrt{2}-\ell_{\text{max}}/2>0$ is no
longer satisfied, i.e., the billboard spills into the street. \ To talk of
asymptotics as $r\rightarrow\infty$, we would perhaps wish to alter the
position of the rectangle center.

\section{Circular Deadbolt with Keyhole}

The initial setting here is the $xy$-plane. \ Let the observer be at
coordinate $\theta$ on the unit circle, centered at the origin. \ Let the
object be the subinterval $[-1/2,1/2]$ of the $y$-axis. \ By the Law of
Cosines, the angle $\Omega(\theta)$ subtended at the observer's eye by the
line segment is
\begin{align*}
\Omega &  =\arccos\left(  \frac{-1+\left[  \cos(\theta)^{2}+\left(
\sin(\theta)-\frac{1}{2}\right)  ^{2}\right]  +\left[  \cos(\theta
)^{2}+\left(  \sin(\theta)+\frac{1}{2}\right)  ^{2}\right]  }{2\sqrt
{\cos(\theta)^{2}+\left(  \sin(\theta)-\frac{1}{2}\right)  ^{2}}\sqrt
{\cos(\theta)^{2}+\left(  \sin(\theta)+\frac{1}{2}\right)  ^{2}}}\right) \\
&  =\arccos\left(  \frac{\frac{3}{4}}{\sqrt{\frac{5}{4}-\sin(\theta)}%
\sqrt{\frac{5}{4}+\sin(\theta)}}\right)  =\arccos\left(  \frac{3}%
{\sqrt{25-16\sin(\theta)^{2}}}\right) \\
&  =\arctan\left(  \frac{4\sqrt{1-\sin(\theta)^{2}}}{3}\right)  =\arctan
\left(  \frac{4}{3}\left\vert \cos(\theta)\right\vert \right)  .
\end{align*}
If $\theta\in\lbrack0,2\pi)$ is uniformly distributed on the circle, then the
first two moments of $\Omega$ are%
\[
\mathbb{E}\left(  \Omega\right)  =\frac{1}{3\pi}\left[  \pi^{2}-6\ln
(2)^{2}-3\operatorname{Li}_{2}\left(  \frac{1}{4}\right)  \right]
=0.6561351817594581454390278...,
\]%
\[
\mathbb{E}\left(  \Omega^{2}\right)  =\operatorname{Li}_{2}\left(  \frac{1}%
{4}\right)  -\operatorname{Li}_{2}\left(  -\frac{1}{4}\right)
=0.5035529367689960607402666...
\]
where
\[
\operatorname*{Li}\nolimits_{2}(\xi)=%
{\displaystyle\sum\limits_{k=1}^{\infty}}
\dfrac{\xi^{k}}{k^{2}}=-%
{\displaystyle\int\limits_{0}^{\xi}}
\dfrac{\ln(1-t)}{t}dt
\]
is the dilogarithm function \cite{Fi-rail}. Also, the probability density
function of $\Omega(\theta)$ is%
\[%
\begin{array}
[c]{ccc}%
\dfrac{6}{\pi}\,\dfrac{1}{\cos(\omega)^{2}\sqrt{16-9\tan(\omega)^{2}}}, &  &
0<\omega<2\arctan\left(  \dfrac{1}{2}\right)
\end{array}
\]
via standard techniques \cite{Pa-rail}.

The final setting here is $xyz$ space. \ Let the observer be at coordinates
$(\theta,\varphi)$ on the unit sphere, centered at the origin. \ Let the
object be the subinterval $[-1/2,1/2]$ of the $y$-axis. \ Note that%
\begin{align*}
\cos(\theta)^{2}\sin(\varphi)^{2}+\left(  \sin(\theta)\sin(\varphi)\mp
\tfrac{1}{2}\right)  ^{2}+\cos(\varphi)^{2}  &  =\sin(\varphi)^{2}\mp
\sin(\theta)\sin(\varphi)+\tfrac{1}{4}+\cos(\varphi)^{2}\\
&  =\tfrac{5}{4}\mp\sin(\theta)\sin(\varphi).
\end{align*}
By the Law of Cosines, the angle $\Omega(\theta,\varphi)$ subtended at the
observer's eye by the line segment is
\begin{align*}
\Omega &  =\arccos\left(  \frac{\frac{3}{4}}{\sqrt{\frac{5}{4}-\sin
(\theta)\sin(\varphi)}\sqrt{\frac{5}{4}+\sin(\theta)\sin(\varphi)}}\right) \\
&  =\arccos\left(  \frac{3}{\sqrt{25-16\sin(\theta)^{2}\sin(\varphi)^{2}}%
}\right)  =\arctan\left(  \frac{4\sqrt{1-\sin(\theta)^{2}\sin(\varphi)^{2}}%
}{3}\right)  .
\end{align*}
If $(\theta,\varphi)\in\lbrack0,2\pi)\times\lbrack0,\pi]$ is uniformly
distributed on the sphere, then the first two moments of $\Omega$ are%
\[
\mathbb{E}\left(  \Omega\right)  =\frac{\pi}{4}%
=0.7853981633974483096156608...,
\]%
\[
\mathbb{E}\left(  \Omega^{2}\right)  =%
{\displaystyle\int\limits_{0}^{1}}
\frac{\eta\arctan\left[  (4/3)\eta\right]  ^{2}}{\sqrt{1-\eta^{2}}}%
d\eta=0.6472381347206737507335484...
\]
and the probability density function of $\Omega(\theta)$ is%
\[%
\begin{array}
[c]{ccc}%
\dfrac{9}{4}\,\dfrac{\tan(\omega)}{\cos(\omega)^{2}\sqrt{16-9\tan(\omega)^{2}%
}}, &  & 0<\omega<\arctan\left(  \dfrac{4}{3}\right)  .
\end{array}
\]
A\ closed-form expression for $\mathbb{E}\left(  \Omega^{2}\right)  $ remains
open. \ Substituting the line segment by a disk of diameter $1$ in the
$xy$-plane, centered at the origin, also gives an unsolved problem.

\section{Sphere around Framing Square}

So far, we have addressed apparent size of lengths (of line segments) and
areas (of planar regions). \ We conclude with apparent magnitude of angles (at
the intersection of two lines). \ 

Again, the setting is $xyz$ space. \ Let the observer be at coordinates
$(\theta,\varphi)$ on the unit sphere, centered at the origin. \ Call this
point $A$. \ Let $B=(1,0,0)$, $C=(0,1,0)$ and $D=(0,0,0)$, that is, a fixed
angle of $\pi/2$ in the base plane. \ The apparent magnitude of $\pi/2$
relative to $A$ is the (dihedral) angle $\alpha$ between normal vectors
$A\times B$, $A\times C$ to the triangular faces $ADB$, $ADC$ respectively:
\[
\alpha=\arccos\left(  \frac{(A\times B)\cdot(A\times C)}{\left\Vert A\times
B\right\Vert \,\left\Vert A\times C\right\Vert }\right)  .
\]
This formula is useful for simulation purposes. \ It's best, however, to
employ the spherical triangle $ABC$ and to recognize that the side $a$
opposite angle $\alpha$ is the constant $\pi/2$. \ From section 1.3 of
\cite{Fn-rail}, the conditional density for angle $\alpha$, given $a=\pi/2$,
is
\[%
\begin{array}
[c]{ccc}%
-\dfrac{1}{\pi}\dfrac{\ln|\cos(\alpha)|}{\sin(\alpha)^{2}}, &  & 0<\alpha<\pi.
\end{array}
\]
\ In particular, a singularity exists at $\alpha=\pi/2$ and
\[%
\begin{array}
[c]{ccc}%
\operatorname*{E}\left(  \alpha\left\vert a=\dfrac{\pi}{2}\right.  \right)
=\dfrac{\pi}{2}, &  & \operatorname*{E}\left(  \alpha^{2}\left\vert
a=\dfrac{\pi}{2}\right.  \right)  =\dfrac{\pi^{2}}{4}+\ln(2)^{2}=\dfrac
{\pi^{2}}{4}+0.480453....
\end{array}
\]

Let us change the vector $C$ from $(1,0,0)$ to $(1/2,\sqrt{3}/2,0)$. \ The
side $a$ opposite angle $\alpha$ is now the constant $\pi/3$. \ It is known
that \cite{Ma-rail}%
\[
\operatorname*{E}\left(  \alpha\left\vert a=\dfrac{\pi}{3}\right.  \right)
=\dfrac{\pi}{3}%
\]
but exact evaluation of%
\begin{align*}
\operatorname*{E}\left(  \alpha^{2}\left\vert a=\dfrac{\pi}{3}\right.
\right)   &  =\frac{\pi^{2}}{12}+%
{\displaystyle\int\limits_{0}^{\pi}}
{\displaystyle\int\limits_{0}^{\pi}}
\arctan\left[  \frac{\left(  1/\sqrt{3}\right)  \sin(y)+\cos(x)\cos(y)}%
{\sin(x)}\right]  ^{2}\frac{\sin(y)}{2\pi}dx\,dy\\
&  =\frac{\pi^{2}}{9}+0.398812...
\end{align*}
appears to be difficult.

Interested readers should examine \cite{Me-rail} for an alternative derivation
of the conditional mean $\operatorname*{E}\left(  \alpha\left\vert a\right.
\right)  $, as well as a higher dimensional analog involving tetrahedra on
$\mathbb{S}^{3}$ (rather than triangles on $\mathbb{S}^{2}$). \ No formula for
the corresponding mean square $\operatorname*{E}\left(  \alpha^{2}\left\vert
a\right.  \right)  $ is known here.

An expression for the solid angle subtended by a polygon can be used to
approximate solid angles for arbitrary piecewise smooth closed curves
\cite{A1-rail, A2-rail}. \ Our treatment of apparent size is based on Euclid
(a visual cone of rays emanating from the eye); this system of thought is
arguably inconsistent with the Renaissance theory of linear perspective
\cite{A3-rail, A4-rail, A5-rail, A6-rail, A7-rail}. \ We welcome
correspondence on this topic, as an opportunity for deeper understanding.

\section{Addendum}

The \textquotedblleft inverse cube law\textquotedblright\ mentioned at the end
of the first vignette appears in other settings. \ We discuss this not with
regard to spherical projection (definition of solid angle) but instead with
respect to planar projection (image formed within a camera). \ Apparent size
is quantified differently than before.

Consider an observer at $(-1,0,1)$ in $xyz$-space, surveying the integer
lattice $\{(m,n,0):m\geq1\}$ in the half $xy$-plane; according to this model,
the observer does not see a uniform grid, but rather its projection into the
$yz$-plane (Figure 3). \ The formula \cite{Nw-rail}%
\[
(x,y,0)\longmapsto\left(  0,\frac{y}{1+x},\frac{x}{1+x}\right)
\]
allows us to calculate the size of the $k^{\text{th}}$ trapezoid in the
(green)\ center strip. \ We find%
\[
\operatorname*{area}=\frac{3+2k}{(1+k)^{2}(2+k)^{2}}\sim\frac{2}{k^{3}}%
\]
as $k\rightarrow\infty$.

Consider instead an observer at $(-\sqrt{2},0,1)$ in $xyz$-space, surveying
the lattice $\left\{  \left(  \frac{m+n}{\sqrt{2}},\frac{-m+n}{\sqrt{2}%
},0\right)  :m\geq1\text{ and }n\geq1\right\}  $ in the quarter $xy$-plane;
the observer here sees a contrasting $yz$-projection (Figure 4). \ The formula%
\[
(x,y,0)\longmapsto\left(  0,\frac{\sqrt{2}(-x+y)}{2+x+y},\frac{x+y}%
{2+x+y}\right)
\]
gives us the size of the $k^{\text{th}}$ quadrilateral in the (gold)\ right
strip. \ We find here%
\[
\operatorname*{area}=\frac{4\sqrt{2}}{(3+k)(4+k)(5+k)}\sim\frac{4\sqrt{2}%
}{k^{3}}%
\]
as $k\rightarrow\infty$.%
\begin{figure}[ptb]%
\centering
\includegraphics[
height=2.0254in,
width=5.9456in
]%
{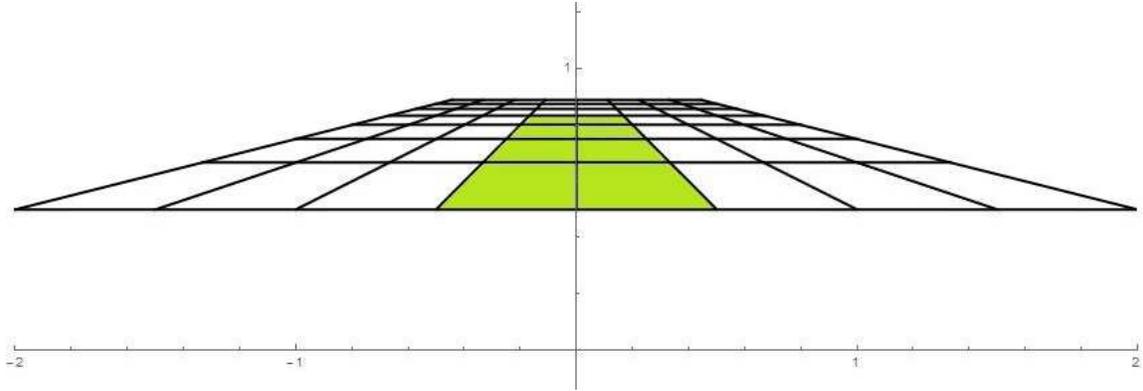}%
\caption{One-point linear perspective: if we drew infinitely many gridsquares,
they would converge at the \textquotedblleft vanishing point\textquotedblright%
\ $(y,z)=(0,1)$.}%
\end{figure}
\begin{figure}[ptb]%
\centering
\includegraphics[
height=3.0234in,
width=5.6628in
]%
{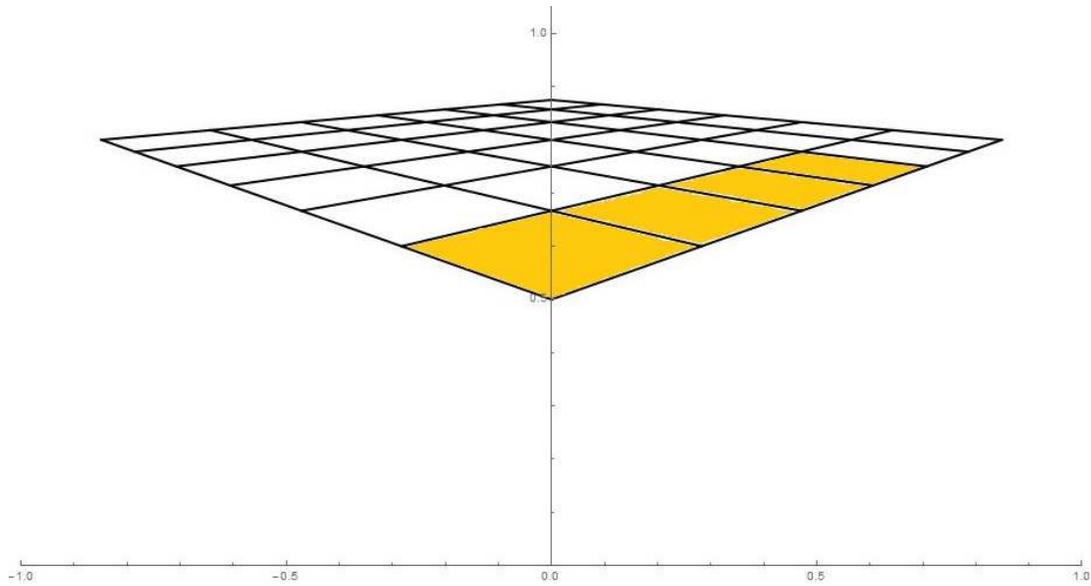}%
\caption{Two-point linear perspective: the \textquotedblleft vanishing
points\textquotedblright\ are at $(y,z)=(\mp\sqrt{2},1)$.}%
\end{figure}

\end{document}